\documentclass[11pt]{amsart}
\usepackage{amsxtra}

\theoremstyle{plain}

 \def\be  {\begin{eqnarray}}
 \def\ee  {\end{eqnarray}}
 \def\ben {\begin{eqnarray*}}
 \def\een {\end{eqnarray*}}

\newcommand \nc {\newcommand}

\newtheorem{theorem}{Theorem}[section]
\newtheorem{lemma}[theorem]{Lemma}
\newtheorem{proposition}[theorem]{Proposition}
\newtheorem{corollary}[theorem]{Corollary}
\newtheorem{definition}[theorem]{Definition}
\newtheorem{example}[theorem]{Example}
\newtheorem{remark}[theorem]{Remark}
\nc \bth[1] { \begin{theorem}\label
{t#1} } \nc \ble[1] {
\begin{lemma}\label{l#1} } \nc \bpr[1] {
\begin{proposition}\label{p#1} } \nc \bco[1] {
\begin{corollary}\label{c#1} } \nc \bde[1] {
\begin{definition}\label{d#1}\rm } \nc \bex[1] {
\begin{example}\label{e#1}\rm } \nc \bre[1] {
\begin{remark}\label{r#1}\rm } \nc \bcon[1] {
\medskip\noindent{\it{Conjecture #1}} } \nc \bqu[1]  {
\medskip\noindent{\it{Question #1}} }
\nc {\ethe} { \end{theorem} }
 \nc {\ele} { \end{lemma} } \nc {\epr}
{ \end{proposition} } \nc {\eco} { \end{corollary} } \nc {\ede} {
\end{definition} } \nc {\eex} { \end{example} } \nc {\ere} {
\end{remark} } \nc {\econ} {\smallskip} \nc {\equ} {\smallskip}
 \nc \thref[1]{Theorem \ref{t#1}}
\nc \leref[1]{Lemma \ref{l#1}} \nc \prref[1]{Proposition
\ref{p#1}} \nc \coref[1]{Corollary \ref{c#1}} \nc
\deref[1]{Definition \ref{d#1}} \nc \exref[1]{Example \ref{e#1}}
\nc \reref[1]{Remark \ref{r#1}}

\def \D {{\mathcal D}}

\def \det { {\mathrm{det}} }

\nc \Gr {Gr} \nc \GRN { \Gr^{(N)} }
\nc \GRA[1] { \Gr_A^{(#1)} }   
\nc \GRAN { \GRA{N} } \nc \GrA[1] { \Gr_A(#1) } \nc \GrAa {
\GrA{\alpha} }
\nc \GRB[1] { \Gr_B^{(#1)} }   
\nc \GRBN { \GRB{N} } \nc \GrB[1] { \Gr_B(#1) } \nc \GrBb {
\GrB{\beta} }
\nc \GRMB[1] { \Gr_{MB}^{(#1)} }   
\nc \GRMBN { \GRMB{N} } \nc \GrMB[1] { \Gr_{MB}(#1) } \nc \GrMBb {
\GrMB{\beta} }

 \def\BB  {\kern+0.1em}
 \def\BBB {\kern+0.15em}
 \def\K   {\kern+0.05em}
 \def\MK  {\kern-0.07em}
 \def\MKK {\kern-0.04em}

\usepackage[hypertex]{hyperref}

\begin{document}

 \vspace{-8ex}

\title[A note on the R. Fuchs's problem for the Painlev\'e equations]
{ A note on the R. Fuchs's problem for the Painlev\'e equations }

\author[Tsvetana Stoyanova]{Tsvetana Stoyanova}

\address{Tsvetana Stoyanova, Department of Mathematics and Informatics,
Sofia University, 5 J. Bourchier Blvd., Sofia 1126, Bulgaria }
\email{cveti@fmi.uni-sofia.bg}

\begin{abstract}
 In this article we consider a first-order completely integrable system of partial
 differential equations $\partial \Phi/\partial x=A(x, t)\,\Phi,\,
 \partial \Phi/\partial t=B(x, t)\,\Phi$ with $\Phi=(\phi_1, \phi_2)^\tau$
 where $A(x, t)$ and $B(x, t)$ are $2 \times 2$ holomorphic matrix functions.
 Under some assumptions we find a variable change by which the system
 $\partial \Phi/\partial x=A(x, t)\,\Phi$ is reduced to an equation independent
 on the variable $t$. As an application we show that the R.Fuchs's
 conjecture for the Painlev\'e equations is true for some algebraic solutions.

    \vspace{0.5ex}

\end{abstract}

 \vspace{-8ex}

\date{}
\maketitle

 \section{ Statement of the result }

  In the first part of this note (see section 2) we consider a first-order linear system
   of partial differential equations
    \be\label{ls1}
      \frac{\partial \Phi(x, t)}{\partial x}
       &=& A(x, t)\,\Phi(x, t)\,,\quad
       A(x, t)=\left(\begin{array}{cc}
           a_{11}(x, t)  &a_{12}(x, t)\\
           a_{21}(x, t)  &a_{21}(x, t)
               \end{array}\right)
       \\[0.4ex]
       \label{ls2}
      \frac{\partial \Phi(x, t)}{\partial t}
       &=& B(x, t)\,\Phi(x, t)\,,\quad
           B(x, t)=\left(\begin{array}{cc}
           b_{11}(x, t)  &b_{12}(x, t)\\
           b_{21}(x, t)  &b_{21}(x, t)
               \end{array}\right)
     \ee
   where $\Phi(x, t)=\left(\begin{array}{c}
            \phi_1(x, t)\\
            \phi_2(x, t)
                    \end{array}
              \right)$
   and $a_{ij}(x, t), b_{ij}(x, t)$ are holomorphic in a domain
   $\D$ in $(x, t)$- space.

   We assume that:
    \begin{description}
    \item[Assumption A.1]\,
     The system \eqref{ls1} - \eqref{ls2} is completely integrable
     in the sense of the Frobenius theorem
      \bth{Fr}(Frobenius)
      The system \eqref{ls1} - \eqref{ls2} is completely integrable
      if and only if
       \be\label{IC}
        \frac{\partial A(x, t)}{\partial t} -
        \frac{\partial B(x, t)}{\partial x} +
        A(x, t)\,B(x, t) - B(x, t)\,A(x, t)=0\,.
       \ee
      \ethe
      The system \eqref{IC} is called the integrability condition
      of the system \eqref{ls1} - \eqref{ls2};
    \item[Assumption A.2]\,
      Both matrices $A(x, t)$ and $B(x, t)$ are traceless;
    \item[Assumption A.3]\,
     All components $a_{ij}(x, t)$ and $b_{ij}(x, t)$ of
     the matrices $A(x, t)$ and $B(x, t)$ respectively are of the form
     $$\,
       \sum F_m(x) G_m(t)\,,
     $$
     where $F_m(x)$ and $G_m(t)$ are holomorphic in the domains
     $\D_1$ and $\D_2$ in $x-$ and $t-$ space respectively,
     such that $\D=\D_1 \times \D_2$.
    \end{description}

   We can transform the system  \eqref{ls1}
   by a standard technique into a second order
   linear equation. It is easy to derive the following classical result:
   \ble{l1}
   Assume that $a_{12}(x, t)$ and $b_{12}(x, t)$ do not vanish
  identically. Then under Assumption A.1 and Assumption A.2
  the first component $\phi_1(x, t)$ of the solution
  $\Phi(x, t)$ of the systems \eqref{ls1} - \eqref{ls2} satisfies
  the following differential equations
   \be\label{s1}
    & &
      \frac{\partial^2 \phi_1(x, t)}{\partial x^2} +
       p_1(x, t)\,\frac{\partial \phi_1(x, t)}{\partial x}
     + q_1(x, t)\,\phi_1(x, t)=0
     \\[1.0ex]\label{s2}
      & &
      \frac{\partial \phi_1(x, t)}{\partial x} =
      p_2(x, t)\,\frac{\partial \phi_1(x, t)}{\partial t} +
      q_2(x, t)\,\phi_1(x, t)\,,
  \ee
  where
   \be\label{coeff}
     p_1(x, t) &=& -\frac{\partial}{\partial x} \log a_{12}(x, t)\,,\\
     q_1(x, t) &=&  \det A(x, t) - \frac{\partial a_{11}(x, t)}{\partial x} +
     a_{11}(x, t)\,\frac{\partial}{\partial x} \log a_{12}(x, t)\,,\nonumber\\
     p_2(x, t) &=& \frac{a_{12}(x, t)}{b_{12}(x, t)}\,,\nonumber\\
     q_2(x, t) &=& a_{11}(x, t) - b_{11}(x, t)\,\frac{a_{12}(x, t)}{b_{12}(x, t)}=
            \nonumber\\[0.2ex]
            &=&
       \frac{1}{2}\left(\frac{\partial}{\partial x} \log b_{12}(x, t)
       -\frac{1}{b_{12}(x, t)}\,\frac{\partial a_{12}(x, t)}{\partial t}
       \right)\,.\nonumber
     \ee
    \ele

    Let the elements $a_{12}(x, t)$ and $b_{12}(x, t)$ have the following
    simple forms
    \ben
       a_{12}(x, t) = g(t)\left[P_1(x) + t\,P_2(x)\right]\,,\quad
       b_{12}(x, t) = g(t)\,P_3(x)\,,
    \een
    where $g(t)$ and $P_i(x)$ are holomorphic functions in
    $\D_2$ and $\D_1$ respectively. Then from the last equation
    of \eqref{coeff} we find that
     \be\label{m0}
      q_2(x, t)=R(x) + M(t)\left[f(x) +  t\,h(x)\right]\,,
     \ee
    where
     \be\label{h-f}
      f(x)=\frac{P_1(x)}{P_3(x)},\quad
      h(x)=\frac{P_2(x)}{P_3(x)}
     \ee
     and $R(x)$ and $M(t)$ are functions dependent only on
     $P_i(x)$ and $g(t)$.

   Our  main result is the following

    \bth{main1}
    Assume that $a_{12}(x, t)$ and $b_{12}(x, t)$ do not vanish
  identically as
    \be\label{m1}
      a_{12}(x, t) = g(t)\left[P_1(x) + t\,P_2(x)\right]\,,\quad
      b_{12}(x, t) = g(t)\,P_3(x)\,,
    \ee
     where $g(t)$ is a holomorphic function in $\D_2$, $P_i(x)$ are
     holomorphic functions in $\D_1$.

   Then under Assumption A.i, $i=1, 2, 3$
       by means of the change the variables
       \be\label{ch}
        \phi_1(x, t) &=& \exp\left(\int R(x) d x\right)\,w(x, t)\,,\\
        \label{ch1}
        \tau &=& t \exp\left(\int h(x)\, d x\right) + \int f(x)
                \exp\left(\int h(x) d x\right) d x\,
              \ee
     with $R(x), h(x)$ and $f(x)$ defined by \eqref{m0}, \eqref{h-f},
     equation \eqref{s1} is reduced to the second order linear differential equation
     \be\label{new1}
       \frac{d^2 w}{d \tau^2} + P(\tau) \frac{d w}{d \tau} + Q(\tau)\,w=0
     \ee
     which is independent on  $t$.
   \ethe

    \bre{r1}
     As we are going to apply \thref{main1} to the R.Fuchs's principle for
     the Painlev\'e equation we formulate it from the point of view
     of the variable $x$. One can rewrite it and \leref{l1} from the point of
     view of the variable $t$. For example: equation \eqref{s1} is considered
     as a second order equation
        \ben
         & &
           \frac{\partial^2 \phi_1(x, t)}{\partial t^2} -
      \frac{\partial}{\partial t} \log b_{21}(x, t)
      \,\frac{\partial \phi_1(x, t)}{\partial t} + \\[0.5ex]
      &+&
      \left[\det B(x, t) + \frac{\partial b_{11}(x, t)}{\partial t}
      - b_{11}(x, t)\,\frac{\partial}{\partial t} \log b_{21}(x, t)
     \right]\,\phi_1(x, t)=0
       \een
       and this equation can be transformed by an appropriate variable change
       into an equation independent on the variable $x$.
    \ere
   The meaning of the transformations \eqref{ch} - \eqref{ch1} is the following.
   Consider the auxiliary system to the quasilinear first-order partial
   differential equation \eqref{s2}
    \ben
     \left|\begin{array}{ccl}
      \dot{x} &= &1\\
      \dot{t} &= &-p_2(x, t)=-t\,h(x) - f(x)\qquad \qquad\qquad\qquad(*)\\
      \dot{\phi}_1 &= &q_2(x, t)\,\phi_1=\left(R(x) + M(t)\,[f(x) + t\,h(x)]\right)\,\phi_1\,,
           \end{array}\right.
    \een
   where $\dot{}\,=\frac{d}{d s}$. Then $\tau$ and $w(x, t)$ defined by
    \ben
      \tau  &=& t \exp\left(\int h(x)\, d x\right) + \int f(x)
                \exp\left(\int h(x) d x\right) d x\,\\
      w(x, t) &=& \phi_1\,\exp\left(- \int R(x)\, d x\right)\,
                  \exp\left(\int M(t)\, d t\right)
    \een
 are two independent first integrals of the system $(*)$. As the exponent
 $\exp\left(\int M(t)\, d t\right)$ does not depend on the variable $x$
 we delete it from the transformation \eqref{ch}.

 In the second part (see section 3) we relate \thref{main1} to the
 R.Fuchs's conjecture for the Painlev\'e equations. The six Painlev\'e
 equations govern the isomonodromic deformations of a linear system
 \eqref{ls1} with rational of $x$ elements $a_{ij}(x, t)$ (the variable $x$
 is usually called the spectral parameter). The general theory of
 isomonodromic deformations ensures that the solution $\Phi(x, t)$ of \eqref{ls1}
 satisfies an additional linear system \eqref{ls2} (the variable $t$
 is called the deformation parameter). The integrability condition \eqref{IC} of
 the systems \eqref{ls1} and \eqref{ls2} leads to the six Painlev\'e
 equations $P_{J}$, \cite{F-N, J-M, Sch}. In such a case we refer to
 \eqref{ls1} - \eqref{ls2} as a linearization of the Painlev\'e equations.
 In \cite{F} R.Fuchs made the following hypothesis: let $y(t)$ is
 an algebraic solution of the Painlev\'e equation, then there exists a
 suitable variable change by which the associated linear equation \eqref{ls1}
 can be transformed into an equation independent on the
 deformation parameter $t$. Moreover in the same paper \cite{F}
 R.Fuchs showed that the linear equation taken at the solutions
 $y=0, 1, \infty, t$ of $P_{IV}$ (obtained as special Picard solutions)
 can be reduced to the hypergeometric equation. The R.Fuchs's conjecture for the Painlev\'e equations was
 recently reinvestigated in \cite{M, O-O}. Utilizing \thref{main1}
 we show that the R.Fuchs's hypothesis is true for some algebraic solutions
 of $P_J$.

   This paper is organized as follows. In the next section we prove
   \thref{main1}. In section 3 we apply \thref{main1} to the R.Fuchs's
   principle for the Painlev\'e equations from the second to the fifth.


   \section{ Proof of the \thref{main1} }

      We will prove the theorem in a general situation when all of the
   functions $M(t), f(x)$ and $h(x)$ are nonconstants. At the end of the proof
   we are going only to list equations \eqref{new1} in the particular cases
   when some of these functions are constants.

   Observe that the assumption of the theorem and the last equation of
   \eqref{coeff} imply
    $$g(t)=B^{-1} \,\exp\left(-2 \int M(t) d t\right)\,,\quad
      R(x)=\frac{1}{2}\left[\frac{P'_3(x)}{P_3(x)} - h(x)\right]$$
    for a  constant $B$.

    Hence
    \ben
      & &
      a_{12}(x, t) = B^{-1} e^{-2 \int M(t) d t}\left[P_1(x) + t\,P_2(x)\right]\,,\,\,
      b_{12}(x, t) = B^{-1} e^{-2 \int M(t) d t}\,P_3(x)\,.
    \een

   On the other hand we have
   \ben
    a_{11}(x, t) - \left[f(x) + t h(x)\right]\,b_{11}(x, t)=
    R(x) + M(t)\left[f(x) + t h(x)\right]\,.
   \een
   Then compatibility condition \eqref{IC} of the linear
   system \eqref{ls1} - \eqref{ls2}
    gives
   \be\label{c1}
     \dot{a}_{11} - b'_{11} + a_{12} b_{21} - b_{12} a_{21}=0
   \ee
   where\,\,\, $\dot{} :=\frac{\partial}{\partial t},\,\,\,
                 ' :=\frac{\partial}{\partial x}$.
   In particular
    \ben
    & &
    a_{21}(x, t) - \left[f(x) + t h(x)\right]\,b_{21}(x, t)=\\
     &=&
     \frac{B\,e^{2 \int M(t) d t}}{P_3(x)}\left[
     (M(t) + b_{11}(x, t))\,h(x) - b'_{11}(x, t) +
     (f(x) + t h(x))\,(\dot{b}_{11}(x, t) + \dot{M}(t))\right]\,.
   \een
    This relation implies
   \ben
    a_{21}(x, t)= \frac{B\,e^{2 \int M(t) d t}}{P_3(x)}\,\widetilde{a}_{21}(x, t)\,,\quad
    b_{21}(x, t)= \frac{B\,e^{2 \int M(t) d t}}{ P_3(x)}\,\widetilde{b}_{21}(x, t)
   \een
  as
   \ben
    & &
    \widetilde{a}_{21}(x, t) - \left[f(x) + t h(x)\right]\,\widetilde{b}_{21}(x, t)
    = \\
    &=&
    (M(t) + b_{11}(x, t))\,h(x) - b'_{11}(x, t) +
     (f(x) + t h(x))\,(\dot{b}_{11}(x, t) + \dot{M}(t))\,.
   \een

   Let us suppose that $b_{11}(x, t)=F_1(x) G_1(t) + F_2(x) G_2(t)$.
   Next, compatibility condition
   \be\label{c3}
    \dot{a}_{21} - b'_{21} + 2 a_{21} b_{11} - 2 a_{11} b_{21}=0
   \ee
   gives equation
     \ben
      & &
     - \frac{\partial \widetilde{b}_{21}(x, t)}{\partial x} +
    \left[f(x) + t h(x)\right]\,\frac{\partial \widetilde{b}_{21}(x, t)}{\partial t} =
     -2 h(x)\,\widetilde{b}_{21}(x, t) -\\[0.2ex]
     &-&
     2 \left[M(t) + F_1(x) G_1(t) + F_2(x) G_2(t)\right]\times\\
     &\times&
     \Big[(f(x) + t h(x))\,\left(\dot{M}(t) + F_1(x)\,\dot{G}_1(t) +
     F_2(x)\,\dot{G}_2(t)\right) +\\
     &+&
     h(x)\,\left(M(t) + F_1(x) G_1(t) + F_2(x) G_2(t)\right) -
     F'_1(x) G_1(t) - F'_2(x) G_2(t)\Big]- \\[0.3ex]
     &-&
     \frac{\partial}{\partial t}
      \Big[(f(x) + t h(x))\,\left(\dot{M}(t) + F_1(x)\,\dot{G}_1(t) +
     F_2(x)\,\dot{G}_2(t)\right) +\\
     &+&
     h(x)\,\left(M(t) + F_1(x) G_1(t) + F_2(x) G_2(t)\right) -
     F'_1(x) G_1(t) - F'_2(x) G_2(t)\Big]\,.
   \een
   We can write the general solution of this quasilinear first-order
   partial differential equation as
   \ben
    \widetilde{b}_{21}(x, t) &=&  e^{ 2 \int h(x) d x}\,F(\tau) -
    \left[F_1(x) G_1(t) + F_2(x) G_2(t) + M(t)\right]^2 -\\
     &-&
    \frac{\partial}{\partial t}\,
    \left[F_1(x) G_1(t) + F_2(x) G_2(t) + M(t)\right]=\\
     &=&
     e^{ 2 \int h(x) d x}\,F(\tau) -
    \left[b_{11}(x, t) + M(t)\right]^2 -
    \frac{\partial}{\partial t}\,
    \left[b_{11}(x, t) + M(t)\right]
    \een
 for an arbitrary holomorphic function $F$ of $\tau$.
 Then the initial equation \eqref{s1} is transformed
 to equation
  \be\label{EQ}
    \frac{d^2 w}{d \tau^2} - F(\tau)\,w=0\,,
  \ee
 which is independent on the variable $t$.

 One can show, in a similar way, that when
 $b_{11}(x, t)= \sum F_i(x)\,G_i(t)$ equation \eqref{s1}
 is again reduced to equation \eqref{EQ}.
 This proves the theorem in a general situation.

 We finish the proof with a list of particular situations.
 If $M(t)\equiv M$ for a nonzero constant $M$ then:
 \begin{itemize}
  \item\,
  if  both $f(x)$ and $h(x)$ do not vanish identically,
  no matter they are constants or not,
  then as above equation \eqref{s1} gets into equation
  \eqref{EQ};
  \item\,
  if $f(x)$ does not vanish identically and $h(x) \equiv 0$
  then equation \eqref{s1} turns into equation
   \be\label{EQ1}
      \frac{d^2 w}{d \tau^2} - \left[ M^2 + F(\tau)\right]\,w=0
   \ee
  for an arbitrary holomorphic function $F$ of
  $\tau=t + \int f(x) d x$;
  \item\,
  if $h(x)$ does not vanish identically and $f(x) \equiv 0$ then
   \ben
    & &
    a_{12}(x, t)= B^{-1} e^{-2 M t}\,t^{A+1}\,P_2(x)\,,\quad
    b_{12}(x, t)= B^{-1} e^{-2 M t}\,t^A\,P_3(x)\\[0.3ex]
    & &
    R(x)=\frac{1}{2}\left[\frac{P'_3(x)}{P_3(x)}
    - (A+1) h(x)\right]
   \een
   for constants $B$ and $A$.
   Equation \eqref{s1} is transformed to equation
      \be\label{EQ2}
      \frac{d^2 w}{d \tau^2} - \frac{A}{\tau}\,\frac{d w}{d \tau} - F(\tau)\,w=0
   \ee
  for an arbitrary holomorphic function $F$ of $\tau=t\exp\left(\int h(x) d x\right)$.
\end{itemize}

  If $f(x)$ does not vanish identically then:
  \begin{itemize}
   \item\,
   if $M(t) \equiv 0, h(x)\equiv 0$ then
    \ben
      & &
     a_{12}(x, t)=B^{-1} e^{A t}\,P_1(x)\,,\quad
     b_{12}(x, t)=B^{-1} e^{A t}\,P_3(x)\\[0.3ex]
     & &
     R(x)=\frac{1}{2}\left[\frac{P'_3(x)}{P_3(x)}
     - A f(x)\right]
    \een
    for  constants $B$ and $A$.
    Equation \eqref{s1} turns into equation
     \be\label{EQ3}
      \frac{d^2 w}{d \tau^2} - A\,\frac{d w}{d \tau} - F(\tau)\,w=0
    \ee
  for an arbitrary holomorphic function $F$ of $\tau=t + \int f(x) d x$;
    \item\,
    We make note that the situation when both $f(x)$ and $h(x)$
    do non vanish identically but $M(t) \equiv 0$ is impossible.
  \end{itemize}

   If $h(x)$ does not vanish identically then:
   \begin{itemize}
    \item\,
    if $M(t) \equiv 0, f(x) \equiv 0$ then
    \ben
      & &
     a_{12}(x, t)=B^{-1} t^{A+1}\,P_2(x)\,,\quad
     b_{12}(x, t)=B^{-1} t^A\,P_3(x)\\[0.3ex]
     & &
      R(x)=\frac{1}{2}\left[\frac{P'_3(x)}{P_3(x)}
    - (A+1) h(x)\right]
    \een
    for  constants $B$ and $A$.
    Equation \eqref{s1} is transformed to equation \eqref{EQ2}
    for an arbitrary holomorphic function $F$ of $\tau=t\exp\left(\int h(x) d x\right)$.
   \end{itemize}
   This proves the theorem. \qed

  We end this section in similar to \leref{l1} and \thref{main1}
  results about the second component $\phi_2(x, t)$ of the solution
  $\Phi(x, t)$ of the systems \eqref{ls1} - \eqref{ls2}.

      \ble{l2}
   Assume that $a_{21}(x, t)$ and $b_{21}(x, t)$ do not vanish
  identically. Then under Assumption A.1 and Assumption A.2
  the second component $\phi_2(x, t)$ of the solution
  $\Phi(x, t)$ of the systems \eqref{ls1} - \eqref{ls2} satisfies
  the following differential equations
   \be\label{s1-2}
    & &
      \frac{\partial^2 \phi_2(x, t)}{\partial x^2} -
      \frac{\partial}{\partial x} \log a_{21}(x, t)
      \,\frac{\partial \phi_2(x, t)}{\partial x} + \\[0.5ex]
      &+&
      \left[\det A(x, t) + \frac{\partial a_{11}(x, t)}{\partial x}
      - a_{11}(x, t)\,\frac{\partial}{\partial x} \log a_{21}(x, t)
     \right]\,\phi_2(x, t)=0\nonumber
     \\[1.0ex]\label{s2-2}
      & &
      \frac{\partial \phi_2(x, t)}{\partial x} =
      \frac{a_{21}(x, t)}{b_{21}(x, t)}\,\frac{\partial \phi_2(x, t)}{\partial t} -
      \left[a_{11}(x, t) - b_{11}(x, t)\,\frac{a_{21}(x, t)}{b_{21}(x, t)}\right]\,\phi_2(x, t)\,.
  \ee
   \ele

     \bth{main2}
    Assume that $a_{21}(x, t)$ and $b_{21}(x, t)$ do not vanish
  identically as
    \be\label{m2}
      & &
     a_{21}(x, t) = g(t)\left[P_1(x) + t\,P_2(x)\right]\,,\quad
     b_{21}(x, t) = g(t)\,P_3(x)\,,\\[0.5ex]
     & &
     a_{11}(x, t) - b_{11}(x, t)\,\frac{a_{21}(x, t)}{b_{21}(x, t)} =
     R(x) + M(t)\left[f(x) +  t\,h(x)\right]\nonumber
    \ee
     where $M(t)$ is a function of $t$, $P_i(x)$ are functions of $x$
     and $h(x)=\frac{P_2(x)}{P_3(x)},\,
          f(x)=\frac{P_1(x)}{P_3(x)}$.

      Then under Assumption A.i, $i=1, 2, 3$
       by means of the change the variables
       \be\label{ch-1}
        \phi_2(x, t) &=& \exp\left(-\int R(x) d x\right)\,w(x, t)\,,\\
        \label{ch1-1}
        \tau &=& t \exp\left(\int h(x)\, d x\right) + \int f(x)
                \exp\left(\int h(x) d x\right) d x\,
              \ee
     equation \eqref{s1-2} is reduced to the second order linear differential equation
     \be\label{new1-2}
       \frac{d^2 w}{d \tau^2} + P(\tau) \frac{d w}{d \tau} + Q(\tau)\,w=0
     \ee
     which is independent on  $t$.
   \ethe


 \section{ The R. Fuchs's principle for the Painlev\'e equations }

 In this section applying \thref{main1} and \thref{main2},
 we show that the R.Fuchs's conjecture is true for some algebraic
 solutions of the Painlev\'e equations from the second to the fifth.
 Unfortunately the particular form of
 $q_2(x, t)=a_{11}(x, t) - b_{11}(x, t)\,a_{12}(x, t)/b_{12}(x, t)$
 in \eqref{coeff} restricts our applications very much. On the other hand
 as the Assumption A.3 is not fulfilled for the sixth Painlev\'e
 equation we are not going to consider this Painlev\'e equation here.
 In fact we make no claim to try all possible applications of \thref{main1}
 and \thref{main2} in the R.Fuchs's principle for the Painlev\'e
 equations. We just give some examples.


   \subsection{ The R. Fuchs's principle for  $P_{II}$ }
  \subsubsection{ Miwa - Jimbo's linearization }
   Miwa - Jimbo's isomonodromic deformation equations for the
   second Painlev\'e equation \cite{J-M} are
   \be\label{IE-0}
    \frac{\partial \Phi(x, t)}{\partial x}
      &=& A(x, t)\,\Phi(x, t)\,,\\
    A(x, t) &=& \left(
                \begin{array}{cc}
                 1  &0\\
                 0  &-1
                \end{array}
                \right)\,x^2 +
                \left(
                \begin{array}{cc}
                0   &u\\
                -2 u^{-1} z  &0
                \end{array}
                \right)\,x +\nonumber\\[0.3ex]
             &+&
                \left(
                \begin{array}{cc}
                 z + t/2  & -u y\\
                 -2 u^{-1} (\theta + y z)  &-z-t/2
                \end{array}
                \right)\,,\nonumber\\[0.8ex]\label{IE1-0}
   \frac{\partial \Phi(x, t)}{\partial t}
      &=& B(x, t)\,\Phi(x, t)\,,\\
      B(x, t)  &=& \frac{1}{2}\,\left(
                   \begin{array}{cc}
                   1  &0\\
                   0  &-1
                   \end{array}
                   \right)\,x + \frac{1}{2}\,
                   \left(
                   \begin{array}{cc}
                   0   &u\\
                   -2 u^{-1} z  &0
                   \end{array}
                   \right)\,,\nonumber
  \ee
  which $\Phi=\left(
           \begin{array}{c}
              \phi_1\\
              \phi_2
                  \end{array}
             \right)$
  where $y, z$ and $u$ are functions of $t$ and $\theta$ is a
  parameter.

  Integrability condition \eqref{IC} of the linear system
  \eqref{IE-0} - \eqref{IE1-0}
   gives
   \ben
    \frac{d y}{d t}=z + y^2 + \frac{t}{2}\,,\quad
    \frac{d z}{d t}=-2 y z - \theta\,,\quad
   \frac{d}{d t} \log u = -y.
   \een
 Eliminating $z$, we obtain the second Painlev\'e equation $P_{II}$
  \ben
   \frac{d^2 y}{d t^2} = 2 y^3 + t y + \alpha
  \een
 with $\alpha=\frac{1}{2}-\theta$.

 \subsubsection{ R. Fuchs' conjecture for the  solution $y \equiv 0$ }

    Equations \eqref{s1} and \eqref{s2} taken at the solution
   $y \equiv 0,\, \theta=1/2\, (\alpha=0)$ of $P_{II}$ are
   \be\label{eq1}
     & &
    \frac{\partial^2 \phi_1}{\partial x^2} -
    \frac{1}{x}\,\frac{\partial \phi_1}{\partial x} -
    x^2\,(x^2 + t)\,\phi_1=0\,,\\[0.3ex]
    & &
    \frac{\partial \phi_1}{\partial x} = 2 x\,
    \frac{\partial \phi_1}{\partial t}\,.\nonumber
   \ee
   Under \thref{main1} we have:
   \ben
     M(t) \equiv 0,\,\,\,f(x)=2 x,\,\,\,
     h(x) \equiv 0,\,\,\,A=0\,,\,\,\,
     R(x) \equiv 0\,.
   \een
   By means of the change of the variable
    \ben
      x^2 + t = \tau
    \een
   equation \eqref{eq1} is reduced to
   equation \eqref{EQ3} taken at $A=0$
    \ben
      \frac{d^2 \phi_1}{d \tau^2} = \frac{\tau}{4}\,\phi_1\,.
    \een
    which is independent on the deformation parameter $t$.
    The last equation after transforma-tion
    $\tau=4^{1/3}\,\xi$ is converted into the Airy
    equation \cite{A-S}
      \be\label{Airy}
      \frac{d^2 \phi_1}{d \xi^2} = \xi\,\phi_1\,.
    \ee
 We notice that when $\theta=1/2$ we have $u=B^{-1}, z=-t/2$ for a constant $B$.

  \subsubsection{ R. Fuchs' conjecture for the  solution $y=-\frac{1}{t}$ }

   Equations \eqref{s1-2} and \eqref{s2-2} taken at the
  solution $y=-\frac{1}{t}, \theta=-1/2\, (\alpha=1)$ of $P_{II}$ are
  are nothing but the equations \eqref{eq1}
   \ben
     & &
    \frac{\partial^2 \phi_2}{\partial x^2} -
    \frac{1}{x}\,\frac{\partial \phi_2}{\partial x} -
    x^2\,(x^2 + t)\,\phi_2=0\,,\\[03.ex]
    & &
    \frac{\partial \phi_2}{\partial x} = 2 x\,
    \frac{\partial \phi_2}{\partial t}\,.
   \een
 which is reduced to the Airy equation \eqref{Airy}, \cite{A-S}
   \ben
      \frac{d^2 \phi_2}{d \xi^2} = \xi\,\phi_2\,.
    \een
 We note that when $\theta=-1/2$ we have $u=B^{-1} t, z=-\frac{t}{2}$ for a constant $B$.


   \subsection{ The R. Fuchs's principle for  $P_{III}$ }

  \subsubsection{ Miwa - Jimbo's linearization }

  Miwa -  Jimbo's isomonodromic deformation  equations for the
  third Painlev\'e equation  are \cite{J-M}
     \be\label{IE}\,\,\,\,\,\,\,\,
    \frac{\partial \Phi(x, t)}{\partial x}
      &=& A(x, t)\,\Phi(x, t)\,,\\
    A(x, t) &=& \frac{1}{2}\,\left(
                \begin{array}{cc}
                 t  &0\\
                 0  &-t
                \end{array}
                \right) + \nonumber\\[0.4ex]
                &+&
                \frac{1}{x}\,\left(
                \begin{array}{cc}
                -\theta_{\infty}/2   &- y w z\\
                -w^{-1}\left((z-t) y + \frac{\theta_{\infty}+\theta_0}{2} \frac{z - t}{z} +\frac{\theta_{\infty}-\theta_0}{2}
                       \right)                &\theta_{\infty}/2
                \end{array}
                \right)+\nonumber\\[0.6ex]
                &+&
                \frac{1}{x^2}\,\left(
                  \begin{array}{cc}
                   z - t/2   &-w z\\
                  w^{-1} (z-t)  &-z + t/2
                 \end{array}\right)\,,
               \nonumber\\[1.8ex]\label{IE1}
                \,\,\,\,\,\,\,\,\,
   \frac{\partial \Phi(x, t)}{\partial t}
      &=& B(x, t)\,\Phi(x, t)\,,\\
      B(x, t)  &=& \frac{1}{2}\,\left(
                \begin{array}{cc}
                 1  &0\\
                 0  &-1
                \end{array}
                \right)\,x + \nonumber\\[0.5ex]
                &+&
                \frac{1}{t}\,\left(
                \begin{array}{cc}
                0   &-y w z\\
                 -w^{-1}\left((z-t) y + \frac{\theta_{\infty}+\theta_0}{2} \frac{z - t}{z} +\frac{\theta_{\infty}-\theta_0}{2}
                       \right)                  &0
                \end{array}
                \right) + \nonumber\\[0.6ex]
                &+&
                \frac{1}{x t}\,\left(
                 \begin{array}{cc}
                  -z+t/2  &w z\\
                  -w^{-1}(z-t)  &z-t/2
                 \end{array}
                 \right)\,,\nonumber
   \ee
    with $\Phi=\left(\begin{array}{c}
                  \phi_1\\
                  \phi_2
                    \end{array}\right)$
   where $y, z$ and $w$ are functions of $t$ and
   $\theta_{\infty}, \theta_0$ are parameters.
   Integrabi-lity condition \eqref{IC} of the
  linear system \eqref{IE} - \eqref{IE1}  gives
   \ben
     t\,\frac{d y}{d t} &=& 4 z y^2 - 2 t y^2 +
     (2 \theta_{\infty} - 1)\,y + 2 t\\
     t\,\frac{d z}{d t} &=& -4 y z^2 + (4 t y - 2 \theta_{\infty}
     +1)\,z + (\theta_0 + \theta_{\infty})\,t\,,\\
     t \frac{d}{d t}\,\log w &=& -
     \frac{(\theta_0 + \theta_{\infty}) t}{z} -2 t y +
     \theta_{\infty}\,.
   \een
  Eliminating $z$ we obtain the third Painlev\'e equation $P_{III}$
  \ben
   \frac{d^2 y}{d t^2}=\frac{1}{y}\,\left(\frac{d y}{d t}\right)^2
   -\frac{1}{t}\,\frac{d y}{d t} +
   \frac{1}{t}\,\left(\alpha\,y^2 + \beta\right) + \gamma\,y^3
   + \frac{\delta}{y}
  \een
  with
  \ben
   \alpha=4 \theta_0,\quad \beta=4 ( 1-\theta_{\infty}),
   \quad \gamma=4, \quad \delta=-4\,.
  \een

 \vspace{0.5cm}

  \subsubsection{ R. Fuchs's principle for the solution $y \equiv 1$ }

  Equations \eqref{s1} and \eqref{s2} taken at the solution
  $y \equiv 1,\,\alpha + \beta=0$\,
   (resp. $\theta_0=\theta_{\infty}-1$) of $P_{III}$ are
   \be\label{ode}
     & &
    \frac{\partial^2 \phi_1}{\partial x^2} +
    \left[\frac{2}{x} - \frac{1}{x+1}\right]\,
    \frac{\partial \phi_1}{\partial x} +\\[0.5ex]
    &+&
    \Big[-\frac{t^2}{4} + \frac{1}{4 (x+1)} +
    \frac{2 t (\theta_{\infty} -1) - 1}{4 x} +
    \frac{8 t^2 +16 (\theta_{\infty}-1) t+3}{16 x^2}
    +\nonumber\\[0.5ex]
     &+&
    \frac{(\theta_{\infty}-1) t}{2 x^3} -
    \frac{t^2}{4 x^4}\Big]\,\phi_1=0\,,\nonumber\\[1ex]
    \frac{\partial \phi_1}{\partial x}
    &=& \frac{t (x+1)}{x (x-1)}\,
    \frac{\partial \phi_1}{\partial t} +
    \left[\frac{\theta_{\infty}-1}{2 x} - \frac{2 \theta_{\infty}-1}{2 (x-1)}
    - \frac{t (x+1)}{x (x - 1)}\right]\,\phi_1\,.\nonumber
   \ee

    Accordingly \thref{main1} we have
    \ben
      M(t)\equiv -1\,,\,\,f(x)\equiv 0\,,\,\,
      h(x)=\frac{x+1}{x (x-1)}\,,\,\,
      A=\theta_{\infty}-1\,,\,\,
      R(x)= \frac{\theta_{\infty}-1}{2 x} - \frac{2 \theta_{\infty} -1}{2 (x-1)}\,.
    \een
   By means of the change of the variables
    \ben
      & &
      \phi_1(x, t)=x^{(\theta_{\infty}-1)/2}\,(x-1)^{(1-2\theta_{\infty})/2}
      \,w(x, t)\\[0.3ex]
      & &
      \tau=\frac{(x-1)^2 t}{x}
    \een
   equation \eqref{ode} is converted to equation \eqref{EQ2}
    \ben
       \frac{d^2 w}{d \tau^2} - \frac{\theta_{\infty}-1}{\tau}\,\frac{d w}{d \tau}
       + \left[-\frac{1}{4} + \frac{\theta_{\infty}-1}{2 \tau}
       + \frac{4 \theta^2_{\infty} -1}{16 \tau^2}\right]\,w=0\,,
   \een
   which is independent on the deformation parameter $t$. Moreover
   after the change
   $$\,
     w=\tau^{(\theta_{\infty}-1)/2}\,v
   $$
   the last differential equation is transformed to
   the Whittaker equation \cite{WH}
     \be\label{W}
    \frac{d^2 v}{d \tau^2} - \left(
    \frac{1}{4} - \frac{\kappa}{\tau} +
    \frac{4 \mu^2 -1}{4 \tau^2}\right)\,v=0
   \ee
  with parameters $\kappa=\frac{\theta_{\infty}-1}{2}\,,\,\,
  \mu^2=\frac{1}{16}$.

  We make note that when $\theta_0=\theta_{\infty}-1$ we have
  $z=(1- 2\theta_{\infty})/4, w=B^{-1} e^{2 t}\,t^{\theta_{\infty}}$
  for a constant $B$.


  \subsection{ The R. Fuchs's principle for  $P_{IV}$ }

  \subsubsection{ Miwa - Jimbo's linearization }

  Miwa -  Jimbo's isomonodromic deformation  equations for the
  fourth Painlev\'e equation  are \cite{J-M}
     \be\label{IE-2}\,\,\,\,\,\,\,\,
    \frac{\partial \Phi(x, t)}{\partial x}
      &=& A(x, t)\,\Phi(x, t)\,,\\
    A(x, t) &=& \left(
                \begin{array}{cc}
                 1  &0\\
                 0  &-1
                \end{array}
                \right)\,x + \left(
                \begin{array}{cc}
                t     &u\\
               2(z - \theta_0 - \theta_{\infty})/u        &-t
                \end{array}
                \right)+\nonumber\\[0.4ex]
                &+&
                \frac{1}{x}\,\left(
                  \begin{array}{cc}
                   - z + \theta_0   &- u y/2\\
                  2 z (z - 2 \theta_0)/u y  &z - \theta_0
                 \end{array}\right)\,,
               \nonumber\\[0.8ex]\label{IE1-2}
                \,\,\,\,\,\,\,\,\,
   \frac{\partial \Phi(x, t)}{\partial t}
      &=& B(x, t)\,\Phi(x, t)\,,\\
      B(x, t)  &=& \left(
                \begin{array}{cc}
                 1  &0\\
                 0  &-1
                \end{array}
                \right)\,x + \left(
                \begin{array}{cc}
                0   &u\\
                2(z - \theta_0 - \theta_{\infty})/u   &0
                \end{array}
                \right)\,,\nonumber
   \ee
    with $\Phi=\left(\begin{array}{c}
                  \phi_1\\
                  \phi_2
                    \end{array}\right)$
    where
    $u, z$ and $y$ are functions of $t$ and $\theta_0, \theta_{\infty}$ are
   parameters.
  Integrabi-lity condition \eqref{IC} of the
  linear system \eqref{IE-2}-\eqref{IE1-2} gives
   \ben
     \frac{d y}{d t} &=& - 4 z +  y^2 +
     2 t\,y + 4 \theta_0\\
     \frac{d z}{d t} &=& - \frac{2}{y}\,z^2 +
     \left(- y + \frac{4 \theta_0}{y}\right)\,z +
     (\theta_0 + \theta_{\infty})\,y\,,\\
      \frac{d}{d t}\,\log u &=& - y - 2 t\,.
   \een
  Eliminating $z$ we obtain the fourth Painlev\'e equation $P_{IV}$
   \ben
    \frac{d^2 y}{d t^2}= \frac{1}{2 y}\,
    \left(\frac{d y}{d t}\right)^2 + \frac{3}{2}\,y^3
    + 4 t\,y^2 + 2(t^2 - \alpha)\,y + \frac{\beta}{y}
   \een
  with
  \ben
   \alpha=2 \theta_{\infty} - 1 ,\quad
   \beta= - 8\,\theta^2_0\,.
  \een

  \subsubsection{ R. Fuchs's principle for the solution $y = -2 t$ }

  Equations \eqref{s1} and \eqref{s2} taken at the solutions
  $y = - 2 t,\, \theta_0=\theta_{\infty}=\frac{1}{2}$
   and $y = - 2 t,\, \theta_0=-\theta_{\infty}=-\frac{1}{2}$
   of $P_{IV}$ are the same and they are
   \be\label{ode-2}
     & &
    \frac{\partial^2 \phi_1}{\partial x^2} +
    \frac{t}{x (x + t)}\,\frac{\partial \phi_1}{\partial x} -
    \left[\frac{1}{4 x^2} + (x + t)^2 +
    \frac{1}{2 x (x + t)}\right]\,\phi_1=0\,,\\[1ex]
     & &
    \frac{\partial \phi_1}{\partial x}
    = \left[1 + \frac{t}{x}\right]\,
    \frac{\partial \phi_1}{\partial t} -
    \frac{1}{2 x}\,\phi_1\,.\nonumber
   \ee
   Under \thref{main1} we have
   \ben
    M(t)\equiv 0\,,\,\,\, f(x)=1\,,\,\,\,
    h(x)=\frac{1}{x}\,,\,\,\,A=0\,,\,\,\,
    R(x)=-\frac{1}{2 x}.
   \een
 By means of the change of the variables
    \ben
      & &
      \phi_1(x, t)=x^{-1/2}\,w(x, t)\\[0.3ex]
      & &
      \tau=t\,x + \frac{x^2}{2}
    \een
   equation \eqref{ode-2} is converted to equation
   \eqref{EQ3} taken at $A=0$
   \be\label{E1-2}
    \frac{d^2 w}{d \tau^2} - w=0\,,
   \ee
   which is independent on the deformation parameter $t$.

  We note that when $\theta_0=-\theta_{\infty}=-\frac{1}{2}$ and
  $\theta_0=\theta_{\infty}=\frac{1}{2}$ we
  have $u=B^{-1}$ and $z=0, z=1$ respectively for a constant $B$.

  \subsubsection{ R. Fuchs's principle for the solution $y = - \frac{2}{3}\,t$ }

  Equations \eqref{s1} and \eqref{s2} taken at the solutions
  $y = -\frac{2}{3}\,t,\,
   \theta_{\infty}=\frac{1}{2},\,\,\theta_0=-\frac{1}{6}$ and
   $y = -\frac{2}{3}\,t,\,
   \theta_{\infty}=\frac{1}{2},\,\,\theta_0=\frac{1}{6}$ of $P_{IV}$ are the
   same and they are
   \be\label{ode-3}
     & &
    \frac{\partial^2 \phi_1}{\partial x^2} +
    \frac{t}{x (3 x + t)}\,\frac{\partial \phi_1}{\partial x} -
    \left[\frac{7}{36 x^2} - \frac{t}{6 x^2 (3 x + t)} +
    \frac{(3 x + t)^2\,( 3 x + 4 t)}{27  x }\right]\,\phi_1=0\,,\\[1.5ex]
     & &
    \frac{\partial \phi_1}{\partial x}
    = \left[1 + \frac{t}{3 x}\right]\,
    \frac{\partial \phi_1}{\partial t} +
    \left[-\frac{1}{6 x} + \frac{2 t}{3}\left(1 + \frac{t}{3 x}\right)\right]\,\phi_1\,.\nonumber
   \ee
  Under \thref{main1} we have
  \ben
   M(t)=\frac{2 t}{3}\,,\,\,\,\, f(x)=1\,,\,\,\,
   h(x)=\frac{1}{3 x}\,,\,\,\,R(x)=-\frac{1}{6 x}\,.
  \een
 By the transformation
  \ben
    \phi_1=x^{-1/6}\,w\,,\quad
    \tau=t x^{1/3} +\frac{3}{4}\,x^{4/3}
  \een
  equation \eqref{ode-3} is reduced to equation \eqref{EQ}
   \ben
   \frac{d^2 w}{d \tau^2}=\frac{4 \tau}{3}\,w\,,
   \een
  which is independent on deformation parameter $t$. After the change
  $\tau=\left(\frac{3}{4}\right)^{1/3}\,\xi$ the last equation
  is converted into the Airy equation \eqref{Airy}, \cite{A-S}
   \ben
    \frac{d^2 w}{d \xi^2}= \xi\,w\,.
   \een
  We make note that when $\theta_{\infty}=\frac{1}{2}, \theta_0=-\frac{1}{6}$
  and $\theta_{\infty}=\frac{1}{2}, \theta_0=\frac{1}{6}$
  we have $u=B^{-1}\,e^{-\frac{2 t^2}{3}}$ and
  $z=-\frac{2}{9}\,t^2, z=-\frac{2}{9}\,t^2+\frac{1}{3}$ respectively
  for a constant $B$.


   \subsection{ The R. Fuchs's principle for  $P_V$ }

  \subsubsection{ Miwa - Jimbo's linearization }

  Miwa -  Jimbo's isomonodromic deformation  equations for the
  fifth Painlev\'e equation  are \cite{J-M}
    \be\label{IE-1}\,\,\,\,\,\,\,\,
    \frac{\partial \Phi(x, t)}{\partial x}
       &=&A(x, t)\,\Phi(x, t)\,,\\
    A(x, t) &=& \frac{1}{2}\,\left(
                \begin{array}{cc}
                 t  &0\\
                 0  &-t
                \end{array}
                \right) + \frac{1}{x}\,\left(
                \begin{array}{cc}
                z + \frac{\theta_0}{2}   &-u(z+\theta_0)\\
                 u^{-1} z                &-z - \frac{\theta_0}{2}
                \end{array}
                \right)+\nonumber\\[0.4ex]
                &+&
                \frac{1}{x-1}\,\left(
                  \begin{array}{cc}
                   -z - \frac{\theta_0+\theta_{\infty}}{2}
                   &u y\left(z+ \frac{\theta_0-\theta_1+\theta_{\infty}}{2}\right)\\
                  -\frac{1}{u y}\left(z + \frac{\theta_0+\theta_1+\theta_{\infty}}{2}\right)
                   &z + \frac{\theta_0+\theta_{\infty}}{2}
                 \end{array}\right)\,,\nonumber
         \ee
                \begin{multline}
   \frac{\partial \Phi(x, t)}{\partial t}
      = B(x, t)\,\Phi(x, t)\,,\\
      B(x, t)  = \frac{1}{2}\,\left(
                \begin{array}{cc}
                 1  &0\\
                 0  &-1
                \end{array}
                \right)\,x + \\[0.2ex]
                +
                \frac{1}{t}\,\left(
                \begin{array}{cc}
                0   &-u\left[z+\theta_0 -y\left(z+\frac{\theta_0-\theta_1+\theta_{\infty}}{2}\right)\right]\\
                u^{-1}\left[z-\frac{1}{y}\left(z+ \frac{\theta_0+\theta_1+\theta_{\infty}}{2}\right)\right]   &0
                \end{array}
                \right)\,,\label{IE1-1}
   \end{multline}
     with $\Phi=\left(\begin{array}{c}
                  \phi_1\\
                  \phi_2
                    \end{array}\right)$
  where
    $y, z$ and $u$ are functions of $t$ and $\theta_0, \theta_1$ and
    $\theta_{\infty}$ are parameters.
  Integrability condition \eqref{IC} of the
  linear system \eqref{IE-1}-\eqref{IE1-1} gives
   \ben
     t\,\frac{d y}{d t} &=& t y - 2 z (y-1)^2 -
     (y-1)\left(\frac{\theta_0 -\theta_1+\theta_{\infty}}{2}\, y  -
     \frac{3 \theta_0 + \theta_1+\theta_{\infty}}{2}\right)\,,\\
     t\,\frac{d z}{d t} &=&  y z\left(z + \frac{\theta_0-\theta_1+\theta_{\infty}}{2}\right)
     - \frac{z+\theta_0}{y}\left(z+\frac{\theta_0+\theta_1+\theta_{\infty}}{2}\right)\,,\\
     t \frac{d}{d t}\,\log u &=& - 2 z - \theta_0 +
     y\left(z + \frac{\theta_0 -\theta_1 + \theta_{\infty}}{2}\right)
      + \frac{1}{y}\left(z + \frac{\theta_0+\theta_1+\theta_{\infty}}{2}\right)\,.
   \een
  Eliminating $z$ we obtain the fifth Painlev\'e equation $P_{V}$
  \ben
   \frac{d^2 y}{d t^2}=\left(\frac{1}{2 y} + \frac{1}{y-1}\right)\,
   \left(\frac{d y}{d t}\right)^2 - \frac{1}{t}\,\frac{d y}{d t}
   +\frac{(y-1)^2}{t^2}\,\left(\alpha y + \frac{\beta}{y}\right) +
   \frac{\gamma y}{t} + \frac{\delta y (y+1)}{y-1}
  \een
  with
  \ben
   \alpha=\frac{1}{2}\left(\frac{\theta_0-\theta_1+\theta_{\infty}}{2}\right)^2,\quad
   \beta=-\frac{1}{2}\left(\frac{\theta_0-\theta_1-\theta_{\infty}}{2}\right)^2,\quad
   \gamma=1 - \theta_0 - \theta_1, \quad \delta=-\frac{1}{2}\,.
  \een

 \vspace{0.5cm}

  \subsubsection{ R. Fuchs's principle for the solution  $y =1-\frac{t}{\theta_1-1}$}

  Equations \eqref{s1} and \eqref{s2}
  taken at the
   solution $y =1-\frac{t}{\theta_1-1}, \theta_0=0, \theta_1+\theta_{\infty}=2$
   of $P_V$ are
   \be\label{ode-1}
     & &
    \frac{\partial^2 \phi_1}{\partial x^2} +
    \frac{1}{x-1}\,\frac{\partial \phi_1}{\partial x} -
    \left[\frac{t^2}{4}  +
    \frac{ t (\theta_1-1)}{2 (x-1)} +
    \frac{\theta^2_1}{4 (x-1)^2}\right]\,\phi_1=0\,,\\[1ex]
    \frac{\partial \phi_1}{\partial x}
    &=& \frac{t}{x-1}\,
    \frac{\partial \phi_1}{\partial t} -
    \left[\frac{2 - \theta_1}{2 (x-1)}
    + \frac{t}{2 (x-1)}\right]\,\phi_1\,.\nonumber
   \ee
  Accordingly \thref{main1} we have
  \ben
    M(t)=-\frac{1}{2},\,\,\, f(x) \equiv 0 ,\,\,\,
    h(x)=\frac{1}{x-1},\,\,\, A=1 - \theta_1\,,\,\,\,
    R(x)=\frac{\theta_1-2}{2(x-1)}\,.
  \een
 By means of the change of the variables
    \ben
      & &
      \phi_1(x, t)=(x-1)^{(\theta_1-2)/2} \,w(x, t)\,,\\[0.3ex]
      & &
      \tau=t (x-1)
    \een
   equation \eqref{ode-1} is converted to equation \eqref{EQ2}
   \ben
    \frac{d^2 w}{d \tau^2} - \frac{1 - \theta_1}{\tau}\,\frac{d w}{d \tau} +
     \left[-\frac{1}{4} + \frac{1-\theta_1}{2 \tau}
    + \frac{1 - \theta_1}{\tau^2}\right]\,w=0\,,
   \een
   which is independent on the deformation parameter $t$.
   Moreover, we can apply the transformation
   \ben
    w=\tau^{(1-\theta_1)/2}\,v
   \een
  to the last differential equation and reduce it to the Whittaker equation
  \eqref{W}, \cite{WH},
 with parameters $\kappa=(1-\theta_1)/2,\,\mu^2=\theta^2_1/4$.

 We make note that when $\theta_0=0,\,\theta_1+\theta_{\infty}=2,\,
 y=1-\frac{t}{\theta_1-1}$ we have
 $u=\frac{B^{-1} t^{2-\theta_1}\,e^{-t}}{\theta_1 -1 -t}$ and $z\equiv 0$ for a constant $B$.

  \subsubsection{ R. Fuchs's principle for the solution  $y \equiv -1$}

   In \cite{K-O} Kazuo Kaneko and Yousuke Ohyama show that R. Fuchs's
   conjecture is true for the rational solution
   $y \equiv -1$ for $\theta_0=\theta_1=1/2$ and arbitrary $\theta_{\infty}$
   of the fifth Painlev\'e equation. We shall give this example as an
   application of \thref{main1}. We remark that our transformations
   are slightly different from these in \cite{K-O}.

   Equations \eqref{s1} and \eqref{s2} taken at the solution
   $y \equiv -1$ for $\theta_0=\theta_1=1/2$ and arbitrary $\theta_{\infty}$
   (see also \cite{K-O}) of $P_V$ are
    \be\label{ode1-2}
      & &
     \frac{\partial^2 \phi_1}{\partial x^2} +
     \left[\frac{1}{x} + \frac{1}{x-1} -
     \frac{2 t}{2(1-\theta_{\infty}) + 2 x t - t}\right]\,
     \frac{\partial \phi_1}{\partial x} +\\[0.4ex]
      &+&
     \Big[-\frac{t^2}{4} - \frac{1}{16 x^2} - \frac{1}{16 (x-1)^2} +
     \frac{4 (\theta_{\infty}-1) t^2}{t^2 - 4(1 - \theta_{\infty})}\cdot\,
     \frac{1}{2 (1 - \theta_{\infty}) + 2 t x - t} + \nonumber\\[0.8ex]
      &+&
     \left(\frac{t^2}{16} + \frac{\theta_{\infty} t}{4} - \frac{t}{4}
     + \frac{1 - \theta_{\infty}}{t - 2(1 - \theta_{\infty})}
     + \frac{2 \theta^2_{\infty} - 4 \theta_{\infty} +3}{8}\right)\,
     \frac{1}{x} +\nonumber\\[0.8ex]
      &+&
      \left(-\frac{t^2}{16} + \frac{\theta_{\infty} t}{4} - \frac{t}{4}
      + \frac{1 - \theta_{\infty}}{t + 2 (1 - \theta_{\infty})} -
      \frac{ 2 \theta^2_{\infty} - 4 \theta_{\infty} + 3}{8}\right)\,
      \frac{1}{x-1}\Big]\,\phi_1\,,\nonumber\\[2.5ex]
      & &
      \frac{\partial \phi_1}{\partial x} =
      \left[\frac{1 - \theta_{\infty}}{x (x-1)} + \frac{t}{2}\cdot\,
      \frac{2 x - 1}{x (x-1)}\right]\,\frac{\partial \phi_1}{\partial t}
      + \nonumber\\[0.3ex]
      &+&
      \left[- \frac{2 x-1}{4 x (x-1)} - \frac{1}{4}\,\left(
      \frac{1 - \theta_{\infty}}{x (x-1)} + \frac{t}{2}\cdot\,
      \frac{2 x - 1}{x (x-1)}\right)\right]\,\phi_1\,.\nonumber
    \ee
   Under \thref{main1} we have
   \ben
     M(t) \equiv -\frac{1}{4},\,\,\,\, f(x)=\frac{1-\theta_{\infty}}{x (x-1)},\,\,\,\,
     h(x)=\frac{1}{2}\,\left(\frac{1}{x} + \frac{1}{x-1}\right),\,\,\,
     R(x)=-\frac{1}{4 x} - \frac{1}{4 (x-1)}\,.
   \een
     By means of the change of the variables
    \ben
      & &
      \phi_1(x, t)=\left(x (x-1)\right)^{-1/4} \,w(x, t)\,,\\[0.3ex]
      & &
      \tau=t \left(x (x-1)\right)^{1/2} - (1 - \theta_{\infty})\,
      \log \frac{\sqrt{x} - \sqrt{x-1}}{\sqrt{x} + \sqrt{x-1}}
    \een
   equation \eqref{ode1-2} is converted to equation \eqref{EQ}
   \ben
     \frac{d^2 w}{d \tau^2} - \frac{1}{4}\,w=0
   \een
   which is independent on the deformation parameter $t$.

  We note that when $\theta_0=\theta_1=1/2$ and $y \equiv -1$ we have
  $z=-(t + 2 + 2 \theta_{\infty})/8$ and $u=B^{-1}\,e^{t/2}$  for a constant $B$.

     \subsubsection{ Kitaev's linearization }

  Kitaev's isomonodromic deformation  equations for the
  dege-nerate fifth Painlev\'e equation with $\delta=0$ are \cite{K}
     \be\label{IE-1-1}\,\,\,\,\,\,\,\,
    \frac{\partial \Phi(x, t)}{\partial x}
      &=& A(x, t)\,\Phi(x, t)\,,\\
    A(x, t) &=& \left(
                \begin{array}{cc}
                 0  &0\\
                 t  &0
                \end{array}
                \right) + \frac{1}{x}\,\left(
                \begin{array}{cr}
                a_1   &a_2\\
                a_3   &-a_1
                \end{array}
                \right)+
                \frac{1}{x-1}\,\left(
                  \begin{array}{cr}
                   b_1     &b_2\\
                   b_3     &-b_1
                 \end{array}\right)\,,
               \nonumber\\[0.8ex]\label{IE1-1-1}
                \,\,\,\,\,\,\,\,\,
   \frac{\partial \Phi(x, t)}{\partial t}
      &=& B(x, t)\,\Phi(x, t)\,,\\
      B(x, t)  &=& \left(
                \begin{array}{cc}
                 0  &0\\
                 x  &0
                \end{array}
                \right) +
                \frac{1}{t}\,\left(
                \begin{array}{cc}
                a_1+b_1   &a_2+b_2\\
                a_3+b_3   &-(a_1+b_1)
                \end{array}
                \right)\,,\nonumber
   \ee
 where $a_i, b_i$ are functions of $t$. Let $y=y(t)$ is an
 arbitrary solution of $P_V$ with parameters
  \ben
   \alpha=\frac{\theta^2_0}{8},\quad \beta=-\frac{\theta^2_1}{8},\quad
   \gamma=\theta_{\infty},\quad \delta=0.
  \een
  Then $a_i, b_i$ are defined as follows
   \ben
     & &
     a_2=\frac{\theta_{\infty}}{2 (y-1)}\,,\,\,\,
    \frac{d}{d t}\left(t \frac{d}{d t} \log a_2\right)=
     \theta_{\infty}\, \frac{d}{d t}\left(\frac{a_1}{a_2}\right)
     + 2 a_2 + \theta_{\infty}\,,\\[0.5ex]
     & &
    a_3=\frac{1}{a_2}\,\left(-\frac{\theta^2_0}{16} - a^2_1\right)\,,\,\,\,
    b_1=\frac{t}{2}\,\frac{d}{d t} \log a_2 - a_1\,\left(
    1 + \frac{\theta_{\infty}}{2 a_2}\right)\,,\\[0.5ex]
     & &
     b_2=- a_2 - \frac{\theta_{\infty}}{2}\,,\,\,\,
     b_3=- \frac{t}{a_2}\,\frac{d}{d t} a_1 - a_3\,
     \left(1 + \frac{\theta_{\infty}}{2 a_2}\right)\,.
   \een

    \subsubsection{ R. Fuchs's principle for the solution  $y =1+ \kappa \sqrt{t}$}

  In this subsection we show that R. Fuchs's principle is true for the
  algebraic solution $y=1 + \kappa \sqrt{t},\,\alpha=\mu,\,
  \beta=-1/8,\,\gamma=- \mu \kappa^2,\,\delta=0$ for arbitrary constants $\kappa$
  and $\mu$ of $P_V$.

  To apply \thref{main1} we make the transformation
   \ben
    t \longmapsto z^2\,.
   \een
   Then equations \eqref{s1} and \eqref{s2} taken at the solution
   $y=1 + \kappa z$ are
   \be\label{ode1-1}
     & &
    \frac{\partial^2 \phi_1}{\partial x^2} +
    \left[\frac{1}{x} + \frac{1}{x-1} - \frac{\kappa z}{\kappa z x +1}\right]\,
    \frac{\partial \phi_1}{\partial x} +\\[0.5ex]
     &+&
    \Big[\frac{\mu}{2 x^2} - \frac{1}{16 (x-1)^2} +
    \frac{4 \mu \kappa z + 2 \mu -1}{4 x} +
    \frac{\kappa^2 z^2}{4 (\kappa z + 1) (1 + \kappa z x)} -\nonumber\\[1.3ex]
     &-&
    \frac{2 \mu \kappa^3 z^3 + 6 \mu \kappa^2 z^2 + 6 \mu \kappa z + 2 \mu -1}
         {4 (\kappa z +1) (x-1)}\Big]\,\phi_1=0\nonumber\\[1.5ex]
    & &
    \frac{\partial \phi_1}{\partial x}=\left[
    \frac{1}{2 \kappa x (x-1)} + \frac{z}{2 (x-1)}\right]\,
    \frac{\partial \phi_1}{\partial z} + \nonumber\\[0.5ex]
    &+&
    \left[
    - \frac{1}{4 (x-1)} + \frac{1}{2 z}\,\left(
    \frac{1}{2 \kappa x (x-1)} + \frac{z}{2 (x-1)}\right)\right]\,\phi_1\,.\nonumber
   \ee
   Accordingly \thref{main1} we have
   \ben
    M(s)=\frac{1}{2 z}\,,\,\,\,\, f(x)=\frac{1}{2 \kappa x (x-1)}\,,\,\,\,\,
    h(x)=\frac{1}{2 (x-1)}\,,\,\,\,\,R(x)=-\frac{1}{4 (x-1)}\,.
   \een
    By means of the change of the variables
    \ben
      & &
      \phi_1(x, z)=(x-1)^{-1/4} \,w(x, z)\,,\\[0.3ex]
      & &
      \tau=z (x-1)^{1/2} - \frac{i}{2 \kappa}\,
      \log \frac{\sqrt{x-1} - i}{\sqrt{x-1} + i}
    \een
   equation \eqref{ode1-1} is converted to equation \eqref{EQ}
   \ben
     \frac{d^2 w}{d \tau^2} - 2 \mu \kappa^2\,w=0
   \een
   which is independent on the deformation parameter $t$.

   We make note that when
   $y=1 + \kappa z,\,\alpha=\mu,\,
  \beta=-1/8,\,\gamma=- \mu \kappa^2,\,\delta=0$ we have
   $a_1=\left(-\frac{2 z}{\kappa} - z^2 + B\right)\,a_2$ for a constant
   $B$.

   \section{ Concluding Remarks }
    We shall address the generalization of \thref{main1} and new applications of
    this theorem in the forthcoming papers.

\begin{small}
    
\end{small}

\end{document}